\documentclass[11pt,oneside,fleqn,reqno,titlepage]{article}
\usepackage{amsmath}
\usepackage{amssymb}
\usepackage{amsthm}
\usepackage{natbib} \citeindextrue
\usepackage{comment}
\usepackage{url}
\usepackage{bbm}
\usepackage{bm}
\usepackage{multirow}
\usepackage{relsize}
\usepackage[pdftex]{graphicx}
\usepackage{color}

\newcommand{\doublespacerestore}{\addtolength{\baselineskip}{1.25\baselineskip}}

\newcommand{\singlespace}{\addtolength{\baselineskip}{-.5\baselineskip}}

\newcounter{stepno}

\newcommand{\E}{\mathbb{E}}
\newcommand{\bn}{\begin{eqnarray}}
\newcommand{\en}{\end{eqnarray}}
\newcommand{\bns}{\begin{eqnarray*}}
\newcommand{\ens}{\end{eqnarray*}}
\newcommand{\defvarbegin}{\begin{quotation}\vspace{-15pt}\begin{tabbing}}
\newcommand{\defvarend}  {\end{tabbing}\vspace{-10pt}\end{quotation}}

\newcommand{\bnarray}{\begin{equation}\begin{array}{rcll}}
\newcommand{\enarray}{\end{array}\end{equation}}

\newcommand{\textwrap}{\parbox[t]{5.0in}}

\newcommand{\barr}{\begin{array}}
\newcommand{\earr}{\end{array}}

\newcounter{cnum}

\newcommand{\beginalg}{\setcounter{stepno}{1}
                \begin{list}{\bf Step~\arabic{stepno}}
                         {\usecounter{stepno}\settowidth{\labelwidth}{\bf Step~9m}
                \addtolength{\leftmargin}{2\parindent}}
                }
\newcommand{\eg}{\end{list}}

\def \define{\begin{quote}\begin{itemize}}
\def \enddefine{\end{itemize}\end{quote}}

\newlength{\boxedparwidth} 
  {\begin{center} \begin{tabular}{|@{\hspace{.15in}}c@{\hspace{.15in}}|}
                 \hline \\ \begin{minipage}[t]{\boxedparwidth}
                 \setlength{\parindent}{0.0in}}%
   {\end{minipage} \\ \\ \hline \end{tabular} \end{center}}
\newcounter{example}
\newcommand{\argmax}{{\rm arg}\max}

\newcommand{\argmin}{{\rm arg}\min}

\def \ptilde{{\tilde p}}

\def \xtilde{{\tilde x}}

\def \Ctilde{{\tilde C}}

\def \Rtilde{{\tilde R}}
\def \Stilde{{\tilde S}}

\def \Wtilde{{\tilde W}}

\def \Omegatilde{{\tilde \Omega}}

\def \chat{\hat c}

\def \phat{\hat p}

\def \cbar{\bar c}

\def \hbar{\bar h}

\def \Cbar{\bar C}

\def \Vbar{\bar V}

\def \mubar{\bar \mu}
\def \sigmabar{\bar \sigma}

\def \Dcal{{\cal D}}

\def \Fcal{{\cal F}}

\def \Ical{{\cal I}}

\def \Lcal{{\cal L}}

\def \Xcal{{\cal X}}

\begin{document}

\title{The Parametric Cost Function Approximation: A new approach for multistage stochastic programming}
\author{Warren B. Powell \\ Department of Operations Research and Financial Engineering \\
        Saeed Ghadimi \\ University of Waterloo}
\date{\today}

\maketitle

\clearpage

\begin{abstract}
The most common approaches for solving multistage stochastic programming problems in the research literature have been to either use value functions (``dynamic programming") or scenario trees (``stochastic programming") to approximate the impact of a decision now on the future.  By contrast, common industry practice is to use a deterministic approximation of the future which is easier to understand and solve, but which is criticized for ignoring uncertainty.  We show that a parameterized version of a deterministic optimization model can be an effective way of handling uncertainty without the complexity of either stochastic programming or dynamic programming.  We present the idea of a parameterized deterministic optimization model, and in particular a deterministic lookahead model, as a powerful strategy for many complex stochastic decision problems.  This approach can handle complex, high-dimensional state variables, and avoids the usual approximations associated with scenario trees or value function approximations.  Instead, it introduces the offline challenge of designing and tuning the parameterization.  We illustrate the idea by using a series of application settings, and demonstrate its use in a nonstationary energy storage problem with rolling forecasts.
\end{abstract}

\thispagestyle{empty}

\singlespace

\doublespacerestore

\clearpage
\setcounter{page}{1}
\pagestyle{plain}

\section{Introduction}
There is an extensive history in the academic literature of solving sequential decision problems over time, under uncertainty, using the idea of approximating the future using sampled scenario trees, dating to the original paper by \cite{Da55}. Since this time thousands of papers have been written using this approach for making decisions under uncertainty.

A parallel literature has evolved estimating the value of being in a downstream state using statistical approximations estimated using machine learning that have evolved under names such as approximate (or adaptive) dynamic programming, neuro-dynamic programming, and (more recently) reinforcement learning (although this is an umbrella for an entire family of methods).  In the arena of stochastic linear programs, we can approximate the value of being in a state using Bender's cuts, a method known as stochastic dual dynamic programming.  Other methods have evolved under the heading of robust optimization (when used in the setting of sequential decision problems), which optimizes over an uncertainty set.

All of these methods are computationally demanding which has sharply limited their use in practice.  At the same time, when used in the context of complex operational problems (our primary focus in this paper), all of these methods are solving approximations of lookahead models which means that none of them are offering optimal or even asymptotically optimal policies (an optimal solution of an approximate problem is not an optimal policy).

By contrast, there has been a long history in industry of using deterministic optimization models to make decisions that are then implemented in a stochastic setting.  Grid operators use deterministic forecasts of wind, solar and loads to plan energy generation \citep{wallace2003stochastic}; utilities use rolling estimates to plan the storage of natural gas \citep{Lai2008}; airlines use deterministic estimates of flight times to schedule aircraft and crews \citep{lan2006planning}; and retailers use deterministic estimates of demands and travel times to plan inventories \citep{harrison1999multi}. There is an extensive literature using deterministic lookahead models for dynamic vehicle routing problems \citep{pillac2013}. All of us use the solution of deterministic shortest path problems produced by in-vehicle navigation systems to plan our route through stochastic, dynamic traffic networks.  These models have been widely criticized by the research community for not accounting for uncertainty, but this criticism ignores the creative use of parametric modifications that help these models account for uncertainty.

%We make the case that these previous approaches ignore the true problem that is being solved, which is always stochastic.  The so-called ``deterministic models'' used in industry are almost always parametrically modified deterministic approximations, where the modifications are designed to handle uncertainty.  Both the ``deterministic models'' and the ``stochastic models'' (formulated using the framework of stochastic programming) are examples of lookahead policies to solve a stochastic optimization problem with the goal of finding the best policy which is typically tested using a simulator, but which may be field tested in an online environment (the real world).

While these parameterized policies have been dismissed as industrial heuristics in the stochastic optimization community, we argue that they parallel parametric functions that have been studied in statistics for over 100 years and widely used in practice.  The field of parametric statistics requires that a knowledgeable human choose the structure of a parameterized model, which might be linear or nonlinear in the parameters (this includes neural networks), and then uses data and algorithms to choose the parameters.  Our use of parameterized deterministic models similarly requires a domain expert to design the parameterized optimization problem, after which algorithms have to be used to find the best values of the parameters.  With parametric statistics, we need a training dataset to estimate the parameters.  With a parameterized policy, we need a model of the problem that is captured in a dataset, although it is possible to do online training in the field.

In this paper, we characterize these modified deterministic models as {\it parametric cost function approximations} (CFAs) which puts them into the same category as other parameterized policies that have been used in computer science for solving simpler problems under the umbrella of ``policy search.'' Policy search is the same as stochastic search applied to the context of tuning parameterized policies for sequential decision problems.  The difference is that using a parameterized optimization model for a policy allows us to solve dramatically more complex problems than we can with simple analytical functions that are the focus of the policy search literature (\cite{Sutton2018}, \cite{Robots}, \cite{Levine2013}).  Instead of controlling a robot or computer game, parametric cost function approximations can be used to dispatch thousands of trucks, schedule an airline, manage a network of hydroelectric reservoirs, or plan energy generation for a power grid.

Our paper takes these ideas into the arena of multistage stochastic programming problems, which has been dominated in the research community by methods based on lookahead policies using scenario trees or, for special cases, Benders' decomposition (\cite{birge2011introduction,sen2014multistage,bayraksan2009assessing,zhao2013multi}).  Our use of modified linear programs is new to the policy search literature, where ``policies'' are typically parametric models such as linear models (``affine policies''), structured nonlinear models (such as $(s,S)$ policies for inventories) or neural networks (such as \cite{HanE16}). These more classical ``policies'' are limited to scalar or low-dimensional problems which could not be applied to the domain of high-dimensional resource allocation problems.

There are two dimensions to our approach:
\begin{description}
\item{1)} The design of the parameterized lookahead model which serves as the policy for making decisions.
\item{2)} The optimization of the parameters so that the policy performs as well as policy in expectation.
\end{description}
The process of designing the parameterization involves the same art as the design of any statistical model or parametric policy; it requires exploiting the structure of the problem which would be associated with a particular problem domain.  The optimization of the parameters is a more classical algorithmic challenge which we formulate and solve as a stochastic optimization problem (more precisely, a stochastic search problem).  Both dimensions are challenging, but the end result is a model that is no harder to solve than a classical deterministic lookahead model.  By contrast, policies based on stochastic programming produce stochastic lookahead models that are much harder to solve in the field than a parameterized deterministic model.  In addition, we argue below that a parameterized deterministic model may easily produce higher quality solutions.

We feel that classical stochastic programs based on scenario trees suffer from several limitations:
\begin{itemize}
    \item They are much harder to solve than deterministic lookahead models.
    \item Two-stage stochastic programs, which are the most widely used approach in stochastic programming, require a number of modeling approximations, including:
    \begin{itemize}
        \item The replacement of fully sequential, multistage stochastic decision problems with two-stage approximations.
        \item The need to use a sample of outcomes, and often a fairly small sample.
        \item The inability to model interactions between decisions and exogenous information (relevant in some applications).
    \end{itemize}
    \item There are many problems where the effect of uncertainty on the behavior of the policy is well known.  Yet, two-stage stochastic programs using modest numbers of scenarios do not provide any mechanism to capture this intuition.
\end{itemize}

Methods based on value function approximations, which includes SDDP (\cite{Pereira1991}, \cite{birge2011introduction}, \cite{Kall2009}), and approximate dynamic programming (\cite{PowellADP2011}, \cite{Bertsekas2017}) represent an alternative approach for optimizing over a stochastic lookahead model.  There are many complex problems where ADP is not a viable method for solving the full model, but might be a useful approach for solving an approximate stochastic lookahead model.  We revisit this issue later.

%Our interest is in problems for which these methods would not be appropriate.  Later we provide guidance regarding problems that lend themselves to value function approximations (of any form) versus a full optimization over a planning horizon, which is the focus of this paper.

This paper formalizes the idea that an effective way to solve certain classes of complex stochastic optimization problems is to shift the modeling of uncertainty from an approximate lookahead model to the stochastic base model, typically implemented as a simulator but which might also be the real world. Tuning a model in a stochastic simulator makes it possible to handle arbitrarily complex dynamics, avoiding the many approximations that are standard in stochastic programming. This strategy also means that we transition from methods that are hard to solve in the field, to methods that are (relatively) easy to solve in the field, but which require serious research in the laboratory to design and tune model parameterizations.

%The CFAs make it possible to exploit structural properties using domain knowledge.  For example, we could solve a stochastic shortest path problem minimizing the risk of arriving late to determine when to leave, or we can solve a deterministic approximation and then add in a buffer to leave early.  Airlines solve large deterministic models to schedule aircraft, but include buffers to account for weather delays (which have to be tuned through field experience).  Grid operators solve deterministic approximations of the future, but schedule reserve generation to handle uncertainties such as generator failures or, increasingly, variations in energy from wind and solar.

The parametric CFA makes it possible to incorporate problem structure for handling uncertainty. Some examples include:
\begin{itemize}
\item Supply chains handle uncertainty by introducing buffer stocks.
\item Hospitals can handle uncertainty in blood donations and the demand for blood by maintaining supplies of O-minus blood, which can be used by anyone.
\item People will use the estimated time for the shortest path when driving to a destination, but will then leave early to accommodate the uncertainty in the travel times.
\item Grid operators handle uncertainty in generator failures, as well as uncertainty in energy from wind and solar, by requiring reserve generator capacity.
\end{itemize}
Central to our approach is the ability to manage uncertainty by recognizing effective strategies for responding to unexpected events.  We would argue that this structure is apparent in many settings, especially in complex resource allocation problems. We offer that our approach represents an interesting, and very practical, alternative to classical stochastic programming.  We also argue that this approach is a perfectly valid form of stochastic optimization that may easily outperform methods based on solutions of approximate stochastic lookahead models.

%We illustrate the parametric CFA using a number of settings.  We then use an inventory problem that arises in energy storage to address a widely overlooked issue, which is the handling of time-dependent problems where inventory decisions have to be made in the presence of capacity constraints and rolling forecasts.  The presence of rolling forecasts is quite common in many resource allocation problems, but there is only a small number of papers that recognize the proper way to model rolling forecasts in the context of a sequential decision problem, and these are limited to very simple problems.

Our presentation is organized as follows.  Section \ref{sec:canonicalmodel} provides a canonical model for sequential decision problems and describes two strategies for designing policies: policy search, which searches over classes of policies (and policies within a class) to find the ones that work best over time (section \ref{sec:policysearchpolicies} provides more detail), and policies based on lookahead models which make good decisions by creating the best approximation of the impact of a decision now on the future (section \ref{sec:lookaheadpolicies} provides more detail).  Basic cost function approximations fall within the policy search class, but we will be considering hybrids that combine deterministic lookaheads (from the lookahead class) with parameterizations (from the policy search class).  Section \ref{sec:examplecfas} provides a series of examples of cost function approximations, divided between parameterizations of the objective function, and parameterizations of the constraints.  Then, section \ref{sec:cfaenergystorage} illustrates the idea of a parameterized, deterministic lookahead policy using the setting of a time-dependent energy storage problem using rolling forecasts, a common modeling device which has been largely overlooked in the research literature.  Section \ref{sec:closingremarks} offers some closing remarks.

%The modeling framework and an overview of the different classes of policies are given in Section 2. We then formally introduce the parametric CFA approach in Section 3. Algorithms for optimizing policy parameters together with their convergence properties and some theoretical results about structure of the optimization problem in the CFA approach are presented in Section 4. We then specialize the parametric CFA approach for an energy storage application and present some numerical experiments for solving this problem in Section 5. Finally, we conclude the paper in Section 6.
%%%%%%%%%%%%%%%%%%%%%%%%%%%
%%%%%%%%%%%%%%%%%%%%%%%%%%%
\section{Canonical model and solution strategies}
\label{sec:canonicalmodel}
We begin by writing a sequential decision problem as the sequence
\bns
(S_0,x_0,W_1, \ldots, S_t, x_t, W_{t+1}, \ldots, S_T)
\ens
where $S_t$ is our state at time $t$ which includes physical state variables $R_t$ (this captures physical resources such as inventories or the location on a graph), information $I_t$ (this could be prices, speeds, weather), and beliefs $B_t$ which captures what we know about uncertain quantities and parameters.  Let $C_t(S_t,x_t)$ be the cost we incur given what we know in $S_t$ and our decision $x_t$.  $W_{t+1}$ is information that arrives after we make a decision, which can depend on $S_t$ and/or $x_t$.  We make decisions $x_t$ using a method we call a {\it policy} that we write as $X^\pi(S_t)$.  The state variable evolves according to a transition function
\bns
S_{t+1} = S^M(S_t,x_t = X^\pi(S_t), W_{t+1}),
\ens
where the transition can span updating inventories, tracking the movement of a vehicle, updating prices and weather, or updating estimates or beliefs about uncertain quantities and parameters.

Our goal is to find a policy $\pi$ that solves
\begin{equation}
	\min_{\pi \in \Pi}  \E \left\{ \sum^T_{t=0} C_t(S_t, X^{\pi}_t(S_t)) \: \bigg| \: S_0  \right\}, \label{eq:baseobjective}
\end{equation}
where $S_{t+1} = S^M(S_t,X^\pi_t(S_t),W_{t+1})$ and given an information model for $(S_0, W_1, \ldots, W_T)$.  We will sometimes refer to \eqref{eq:baseobjective} as the base model, which will typically be a simulator but is sometimes the real world.

This canonical model is the foundation for a wide range of sequential decision problems.  \cite{Powell2019} describes two fundamental strategies for designing policies: policy search, where we search over classes of functions for making decisions to find the function that works the best over time, and policies based on lookahead approximations which combine the immediate cost of a decision plus an approximation of costs incurred as a result of the decision.

Each of these two strategies can be applied in two different ways, creating four classes of policies which encompass every possible method for making decisions \citep{Powell2019}:
\begin{itemize}
    \item Policy search strategies:
    \begin{description}
        \item{1)} Policy function approximations (PFAs) - These are analytical functions that map states to decisions.
        \item{2)} Cost function approximations (CFAs) - These are parameterized static or single-period optimization problems that yield a decision.
    \end{description}
    \item Policies based on lookahead approximations:
    \begin{description}
        \item{3)} Policies based on value function approximations (VFAs) - A value function approximation estimates the future cost from the next state we land in after making decision $x_t$ now.
        \item{4)} Policies based on approximations of direct lookahead models (DLAs) - Here we form an approximate model over some horizon, which we solve to make a decision now. Lookahead models can be divided into two broad classes:
        \begin{itemize}
            \item Deterministic lookahead models.
            \item Stochastic lookahead models.
        \end{itemize}
    \end{description}
\end{itemize}
We can also create hybrids, as we will below when we create a parameterized, deterministic lookahead policy that is a hybrid DLA/CFA.

Policy function approximations have received considerable attention in computer science under the banner of ``policy search''  \citep{ng2000pegasus,peshkin2000learning,hu2007evolutionary,mannor2003cross}.  Parameterized policies such as order-up-to inventory ordering policies have been studied since 1960 \citep{Clark1960}, along with a host of other specially structured policies in dynamic programming (\cite{putermanmarkov}, \cite{PowellRLSO}[Chapter 14]). Policy search uses classical stochastic search methods applied to parameterized analytical functions for making decisions, drawing on an extensive literature in derivative-based and derivative-free methods dating to 1951. We review this literature in section \ref{sec:algorithmsforpolicysearch}.

Policies based on value functions and value function approximations have been investigated extensively, building on the foundation of Bellman's equation \citep{putermanmarkov,Bertsekas2017} or Hamilton-Jacobi equations \citep{Kirk2012,stengel1986,Sontag1998,sethi2019,lewis2012}.  Value function approximations have been studied under names such as approximate dynamic programming \citep{PowellADP2011,Bertsekas2017} with applications in trucking \citep{SiDaGe09}, rail \citep{BeChPo2014}, and health \citep{Bartroff2010}; reinforcement learning  \citep{sutton1998reinforcement,Sutton2018} with applications in robotics \citep{SiBaPo04,Robots} and games \citep{Fu2017};  adaptive dynamic programming \citep{LewisVrabie2009,Wang2009}, neuro-dynamic programming \citep{Neuro_DP} and heuristic dynamic programming \citep{SiBaPo04}.  The stochastic programming community has developed the idea of approximating value functions using Bender's cuts under the heading of stochastic dual dynamic programming (SDDP) \citep{Pereira1991,Shapiro2014} with applications in hydroelectric planning \citep{Shapiro2011,Philpott2000}.

Direct lookahead policies represent an umbrella for a number of strategies that have been studied under names such as model predictive control \citep{camacho2013model} and stochastic programming \citep{birge2011introduction,Kall2009,bayraksan2009assessing,zhao2013multi,sen2014multistage}.  This approach has been developed in many books and thousands of papers, with applications that include unit commitment (\cite{jin2011modeling}), hydroelectric planning (\cite{carpentier2015managing}), and transportation (\cite{lium2009study}).

Another form of direct lookahead policy has evolved more recently using robust optimization which replaces scenario trees with uncertainty sets \citep{Ben-Tal2009a,Wiesemann2014}.  \cite{Ben-Tal2005} illustrates robust optimization as a lookahead policy for an inventory problem. See \cite{BertsimasBrown2011} for a review of other applications of robust optimization.

Cost function approximations, however, have been largely overlooked by the academic research community.  CFAs are parameterized optimization models, where the specific parameterization is designed (as would be the case with any parametric model in machine learning) to incorporate intuition into how uncertainty would affect the solution.  The use of parameterized optimization models has been widely used in industry, but in an ad hoc manner.

There is one example of a cost function approximation which has received extensive attention from the research literature.  The problem of finding the best performer out of a discrete set of alternatives (drugs, products, ads) is known as the multiarmed bandit problem.  A widely studied class of policies are known as upper confidence bounding, first introduced by \cite{Lai1985}.  A simple version (introduced by \cite{Ka93}) is given by
\bn
X^{UCB}(S_t|\theta) &=& \argmax_{x} \big(\mubar^n_x + \theta \sigmabar^n_x\big), \label{eq:ucbpolicy}
\en
where $\mubar^n_x$ is the current estimate of the performance of discrete alternative $x\in\{x_1, \ldots, x_M\}$ (e.g. the expected sales from advertising product $x$), and $\sigmabar^n_x$ is the standard deviation of $\mubar^n_x$ (note that $S_t = B_t = (\mubar^n_x, \sigmabar^n_x)$).  There is by now an extensive literature proving various regret bounds on the performance of the policy $X^{UCB}(S_t)$ (see e.g. \cite{Bubeck2012}).  This is particularly important for this paper, since $X^{UCB}(S_t)$ is, first and foremost, a class of parametric cost function approximation given that there is an imbedded ``$\argmax_x$'' within the policy.  The regret bounds that have been derived for UCB policies (of which \eqref{eq:ucbpolicy} is just one example) represents a rare set of provable bounds for the quality of a CFA policy.

PFAs and VFAs tend to be limited to problems that are relatively simple, or which enjoy special structure that can be exploited to estimate the required function (the policy or the value function).  Stochastic lookaheads have attracted considerable attention in the research literature, but relatively little of this work has made its way into practice. For this reason the most common approach used in practice for more complex decision problems is a deterministic approximation, which may be either static or single-period optimization models, or deterministic lookaheads.

In this paper, we want to shine a light on the power of using parameterized deterministic models, particularly for the complex problems that often arise in real applications.  We are not going to argue that this is a panacea that can replace all other methods, but we do feel that it is a powerful and overlooked strategy that belongs alongside widely studied (but rarely used) methods such as stochastic programming or approximate dynamic programming.

We next provide an overview of the policies based on policy search where we cover PFAs and static/single-period CFAs (section \ref{sec:policysearchpolicies}), and policies based on lookahead approximations where we cover policies based on VFAs and deterministic or stochastic DLAs (section \ref{sec:lookaheadpolicies}).  Ultimately we are going to take the idea of policy search that originated with PFAs, and apply it to the idea of parameterized optimization models with special emphasis on parameterized, deterministic lookaheads.  We are going to argue that this is often going to be a more effective strategy in practice for complex sequential decision problems than policies based on VFAs or stochastic DLAs.

\section{Policies based on policy search}
\label{sec:policysearchpolicies}
We begin with the principle of creating parameterized policies (broadly defined), which we divide between policy function approximations and cost function approximations.

\subsection{Policy function approximations (PFAs)}
\label{sec:pfas}
Policy function approximations are analytical functions that map states to decisions.  The functions can be lookup tables, parametric functions (linear or nonlinear), or nonparametric (in particular, locally parametric).  PFAs (using any of a wide range of approximation strategies) have been widely studied in the computer science literature under the umbrella of policy search (see e.g., \cite{Sutton2018}[Chapter 13], \cite{hadjiyiannis2011efficient,lillicrap2015continuous,levine2014learning}).  Although limited to relatively simple problems, the fundamental idea of policy search is quite powerful, and an idea that we are going to exploit.

An example of a PFA is a linear decision rule (also known as an affine policy) which can be written
\bn
X^\pi(S_t|\theta) = \sum_{f\in\Fcal} \theta_f \phi_f(S_t). \label{eq:linearpfa}
\en
where $\phi_f(S_t)$ is a feature drawn from the information in $S_t$ and $\theta_f$ is the coefficient for that feature.

We might wish to use a nonlinear model to choose the price $x^{bid}_t$ to bid for a set of keywords to maximize ad-clicks such as
\bns
X^\pi(S_t|\theta) = \frac{e^{\sum_{f\in\Fcal}\theta_f \phi_f(S_t)}}{1+e^{\sum_{f\in\Fcal}\theta_f \phi_f(S_t)}}.
\ens
More generally, we could represent a policy using a neural network, in which case $\theta$ might have many thousands (or millions) of parameters.

We have to choose the type of function $f\in\Fcal$ (for a neural network, $f$ would specify the network structure, number of layers and nodes per layer). For a given function type $f$ we have to choose $\theta\in\Theta^f$.

The choice of function $f$ is the art of policy search (just as it is with machine learning).  If this were a machine learning problem where we are given a training dataset $(x^n,y^n)_{n=1}^N$, we would be solving
\bn
\min_{(f\in\Fcal, \theta\in\Theta^f)}\sum_{n=1}^N(y^n - f(x^n|\theta))^2.\label{eq:machinelearning}
\en
Policy search for PFAs would be written similarly
\bn
\min_{\pi =(f,\theta) \in (\Fcal,\Theta^f)} \E\left\{\sum_{t=0}^T C_t(S_t,X^\pi_t(S_t|\theta^)) \big| S_0\right\}, \label{eq:objfunpolicies}
\en
where $S_{t+1}=S^M(S_t,x_t=X^\pi(S_t|\theta),W_{t+1}$) and where we are given a model for \linebreak $(S_0,W_1, \ldots, W_T)$.  We note that machine learning \eqref{eq:machinelearning} requires a training dataset while policy search for sequential decision problems in \eqref{eq:objfunpolicies} requires a model of the decision problem (costs, constraints, transition function and exogenous information model).

The comparison with machine learning hints at the limitation of PFAs: they only work for relatively simple decision problems.  For example, the response $y$ is typically scalar, although it might be a low dimensional vector.  We could never use a PFA to, say, schedule machines, dispatch a fleet of trucks, find a shortest path or optimize flows in a supply chain.

\subsection{Cost function approximations (CFAs)}
\label{sec:cfas}

The second class of policy is cost function approximations (CFAs) which are parametrically modified optimization problems.  Eventually we are going to include optimization problems that extend into the future, but for now we limit ourselves to static or single-period optimization models. We emphasize that with the notable exception of upper confidence bounding policies (equation \eqref{eq:ucbpolicy}) for multiarmed bandit problems, CFAs have received virtually no attention in the research literature.

We can create a CFA by modifying either the objective function or the constraints.  For this reason, we begin by defining
\bns
\Cbar^\pi_t(S_t,x_t|\theta) &=& \textwrap{the modified objective function as determined by the policy $\pi$, where $\theta$ represents the tunable parameters,}\\
\Xcal^\pi_t(\theta) &=& \textwrap{the modified set of constraints (that is, the feasible region) determined by policy $\pi$, with tunable parameters $\theta$.}
\ens
We might modify the objective function with a linear correction factor which we could write
\bn
\Cbar^\pi_t(S_t,x_t|\theta) = C(S_t,x_t) + \sum_{f\in\Fcal} \theta_f \phi_f(S_t,x_t). \label{eq:linearobjcorrection}
\en
As an illustration of how constraints can be modified, assume that we start with linear constraints
\bn
A_t x_t &=& b_t, \\
x_t & \leq & u_t,\\
x_t & \geq & 0.
\en
We might modify these using
\bn
A^\pi_t(\theta^a) \xtilde_t &   =  &  \theta^b \odot b_t + \theta^c, \label{eq:constraintcorrection1}\\
x_t                         & \leq &  u_t-\theta^u,\label{eq:constraintcorrection2}\\
x_t                         & \geq &  0 + \theta^\ell. \label{eq:constraintcorrection3}
\en
where $\theta^b \otimes b_t$ is the element by element product of the vector $b_t$ with the similarly dimensioned vector of coefficients $\theta^b$, plus a shift vector $\theta^c$.  We can enter schedule slack by parameterizing the matrix $A^\pi_t(\theta^a)$.  We then reduce the upper bounds $u_t$ by a shift vector $\theta^u$, and possibly raise the lower bounds by $\theta^\ell$.  Our constraints are now parameterized by the (possibly high-dimensional) vector $\theta = (\theta^a, \theta^b, \theta^c, \theta^\ell, \theta^u)$.

A parametric CFA policy can then be written
\bn
X^{CFA}(S_t|\theta) = \argmin_{x_t\in\Xcal^\pi_t(\theta)} \Cbar^\pi(S_t,x_t|\theta). \label{eq:cfa}
\en
We provide examples of CFAs in section \ref{sec:examplecfas}.

\subsection{Algorithms for policy search}
\label{sec:algorithmsforpolicysearch}
PFAs and CFAs both involve tuning a vector $\theta$.  Once we have the structure of the policy $f\in\Fcal$, which is typically chosen by a knowledgeable human guided by intuition, we tune $\theta$ using
\bn
\min_{\theta} \E\left\{\sum_{t=0}^T C_t(S_t,X^\pi_t(S_t|\theta)) \big| S_0\right\}. \label{eq:objfuntheta}
\en
The tools for optimizing the parameters using \eqref{eq:objfuntheta} fall under the broad umbrella of stochastic search which can be approached using both derivative-based algorithms, with a literature that dates to \cite{RoMo51}, and derivative-free algorithms, with a literature that dates to \cite{Box1951}.

The derivative-based stochastic optimization literature is extensive, beginning with the body of research building off of \cite{RoMo51} in the 1950s and later (see, e.g., \cite{Dv56}), initially for unconstrained problems.  A separate literature evolved in the context of constrained stochastic gradient problems (\cite{Sh79}, \cite{Ermoliev1983}).  In addition, the simulation-optimization community (see \cite{fu2015handbook}) has developed powerful tools for taking derivatives of simulations (see \cite{glasserman1991gradient}, \cite{ho1992discrete}, \cite{kushner2003stochastic}, \cite{cao2008stochastic});  a nice tutorial is given in \cite{chau2014simulation}. Much of this literature focuses on derivatives of discrete event simulations, but there is an equally extensive literature on methods based on numerical derivatives such as SPSA (\cite{spall2005introduction,NesSpo17,GhaLan12}).  More recently is work on derivatives of parameterized policies for discrete dynamic programs from the reinforcement learning literature under the umbrella of the policy gradient theorem (\cite{SuMcSi2000}, \cite{Sutton2018}[Chapter 13]).

There is a parallel literature in derivative-free algorithms for stochastic search which is equally extensive.  This literature spans active learning problems \citep{Settles2010}, multiarmed bandit problems \citep{gittins2011}, and optimal learning \citep{PoRy2012}. See \cite{PowellRLSO}[Chapter 7] for an overview of this rich field.

While both PFAs and CFAs require parameter tuning, the characteristics of the tuning problems for PFAs and CFAs tend to be quite different.  It is well known that scaling is a major issue in stochastic search.  Consider the linear decision rule in equation \eqref{eq:linearpfa}.  The scaling of each coefficient $\theta_f$ depends heavily on the characteristics of the feature $\phi_f(S_t)$.  By contrast, the coefficients $\theta$ used in a parametric CFA tend to be scaled by the structure of the deterministic optimization model.  In section \ref{sec:cfaenergystorage}, we demonstrate a parametric CFA for a stochastic inventory control problem where the optimal coefficients are all equal to 1.0 if the forecasts are perfect.  Given imperfect forecasts, the optimal coefficients all appear to be in the interval $[0,2]$.

\section{Policies based on lookahead approximation}
\label{sec:lookaheadpolicies}
%CFAs are widely used in industry for complex problems such as scheduling energy generation or planning supply chains, but they have not been studied formally in the research literature.
    %In special cases, PFAs and CFAs may produce optimal policies, although generally we are looking for the best within a class.

The second strategy for creating policies is to construct policies based on approximations of the downstream impact of a decision $x_t$ made while in state $S_t$.  An optimal policy can be written
\bn
X^*_t(S_t) = \argmin_{x_t \in {\cal X}_t} \left(C_t(S_t,x_t) + \E \left\{ \min_{\pi\in\Pi} \E \left\{\sum_{t'=t+1}^T C_{t'}(S_{t'},X^\pi_{t'}(S_{t'})) \big| S_{t+1}\right\} \big| S_t,x_t \right\}\right)\hspace{-.04in}. \label{eq:optimalpolicylookahead}
\en
Equation \eqref{eq:optimalpolicylookahead} is called a lookahead policy.  Not surprisingly this is computationally intractable for any realistic problem (this includes all the problems that we are interested in).

Just as we divided the policy search strategy into two classes (PFAs and CFAs), there are two classes of policies that we can use to approximate equation \eqref{eq:optimalpolicylookahead}.  These are policies based on value function approximations (VFA policies) and policies based on approximations of the direct lookahead approximations (DLA policies).  We describe these in more detail in sections \ref{sec:vfapolicy} and \ref{sec:dlapolicy} below.

\subsection{Policies based on value function approximations}
\label{sec:vfapolicy}
Equation \eqref{eq:optimalpolicylookahead} is basically Bellman's equation although it is more conventional to write
\bn
X^*_t(S_t) = \argmin_{x_t \in {\cal X}_t} \big(C_t(S_t,x_t) + \E \left\{V_{t+1}(S_{t+1})|S_t,x_t\right\}\big) \label{eq:optimalpolicylookahead2}
\en
where
\bn
V_{t+1}(S_{t+1}) = \min_{\pi\in\Pi} \E \left\{\sum_{t'=t+1}^T C_{t'}(S_{t'},X^\pi_{t'}(S_{t'})) \big| S_{t+1}\right\}. \label{eq:optimalpolicylookahead3}
\en
Alternatively, we can write the expression for value functions recursively using
\bn
V_t(S_t) = \min_{x_t}\big(C_t(S_t,x_t) + E  \left\{V_{t+1}(S_{t+1})|S_t,x_t\right\}\big).\label{eq:optimalpolicylookahead4}
\en
Equation \eqref{eq:optimalpolicylookahead4} is the most common way of writing Bellman's equation, but it is mathematically equivalent to equation \eqref{eq:optimalpolicylookahead} assuming that we can compute $V_t(S_t)$ using either \eqref{eq:optimalpolicylookahead3} or \eqref{eq:optimalpolicylookahead4}, which is never the case for the problems that we are interested in.

When $x_t$ is a vector, it is customary to eliminate the expectation in \eqref{eq:optimalpolicylookahead2} and \eqref{eq:optimalpolicylookahead4} by using the post-decision state variable $S^x_t$ (see \cite{PowellADP2011} and \cite{Shapiro2011}).  We then replace the post-decision value function $V^x_t(S^x_t)$ (which we could never compute) with an approximation $\Vbar^x_t(S^x_t)$, giving us the policy
\bn
X^{VFA}_t(S_t) = \argmin_{x_t}\big(C_t(S_t,x_t) + \Vbar^x_t(S^x_t)\big),\label{eq:optimalpolicylookahead5}
\en
where $\Vbar^x_t(S^x_t)$ might be a linear model, separable piecewise linear functions, or Benders cuts.

While this approach has attracted considerable attention in the literature (see e.g., \cite{powell2004learning,PowellADP2011,bertsekas2011dynamic, Sutton2018}), it is limited to a surprisingly narrow set of problems.  For example, while SDDP has been widely studied in the stochastic programming community, applications are generally limited to fairly simple resource allocation problems such as hydroelectric planning problems \citep{Shapiro2011,Philpott2012} which is known as a single layer resource allocation problem (water is the only resource).  For example, you could never use SDDP for dynamic vehicle routing problems, complex inventory planning problems, or dynamic shortest path problems.

As an indication of how easy it is to ``break'' approximate dynamic programming, ADP (or SDDP) is very effective for solving a blood management problem as long as the surgeries requiring blood having to be completed at a particular point in time (this might be some time during a week).  If we have elective surgeries that can be delayed, creating effective value functions becomes dramatically more difficult.

%Other examples of single layer resource allocation problems arise in the optimizing of fleets of trucks moving loads that must be moved at a single point in time \citep{SiDaGe09}, fleets of locomotives moving trains that also have to be moved at a single point in time \citep{BeChPo2014}.  If we allow the loads of freight or trains to be delayed, the challenge of estimating value function approximations becomes much more complicated.

Even more complex problems include dynamic vehicle routing, where you have to optimize the movement of vehicles and the timing of deliveries, or scheduling machines to complete a series of tasks.  In short, value function approximations are effective when there is structure (such as linearity or convexity) that can be exploited.  Later, we are going to illustrate a CFA in the solution of a time-dependent inventory problem with rolling forecasts.  Rolling forecasts create complex, high-dimensional state variables that are completely intractable using methods based on value function approximations.

\subsection{Policies based on direct lookahead approximations}
\label{sec:dlapolicy}
The most commonly used approach used to solve complex time-dependent problems is to solve an approximate lookahead model on a rolling basis.  Here we first create an approximate sequential decision problem that we are going to use as an approximate model of the future.  We represent this approximate sequential decision problem as
\bns
(\Stilde_{tt},\xtilde_{tt},\Wtilde_{tt}, \ldots, \Stilde_{tt'},\xtilde_{tt'},\Wtilde_{t,t'+1}, \ldots),
\ens
where $\Stilde_{tt'}$, $\xtilde_{tt'}$ and $\Wtilde_{t,t'+1}$ are approximations of $S_{t'}$, $x_{t'}$ and $W_{t'}$ for a decision we are making at time $t$.

While there are a variety of strategies for approximating lookahead models, the two that have received the most attention are:
\begin{itemize}
    \item{Deterministic lookaheads} - This is the approach most widely used in practice, but it has a substantial academic following under the umbrella ``model predictive control.'' Using a deterministic lookahead model reduces  equation \eqref{eq:optimalpolicylookahead} to
    \bn
    X^{DLA}_t(S_t) &=& \argmin_{\xtilde_{tt}} \left(C(\Stilde_{tt},\xtilde_{tt}) + \min_{\xtilde_{t,t+1},\ldots,\xtilde_{t,t+H}} \sum_{t'=t+1}^{t+H} C(\Stilde_{tt'},\xtilde_{tt'})\right) \nonumber\\
               &=& \argmin_{\xtilde_{tt},\xtilde_{t,t+1},\ldots,\xtilde_{t,t+H}}  \sum_{t'=t}^{t+H} C(\Stilde_{tt'},\xtilde_{tt'}). \label{eq:deterministicdla}
    \en
    Equation \eqref{eq:deterministicdla} is so widely used it is known under a number of names including rolling horizon procedure, receding horizon procedure, model predictive control, or deterministic direct lookahead (\cite{sethi1991theory,camacho2013model}, \cite{PowellRLSO}[Chapter 19]).
    \item{Stochastic programming} - First introduced by  \cite{Da55}, the most common approach is to replace the fully sequential decision problem in the future with a two-stage approximation which means our sequence of decisions and information looks like
    \bns
    \big(\xtilde_{tt}, (\Wtilde_{t,t+1}(\omega),\Wtilde_{t,t+2}(\omega), \ldots, \Wtilde_{t,t+H}(\omega)),(\xtilde_{t,t+1}(\omega),\xtilde_{t,t+1}(\omega),\ldots, \xtilde_{t,t+H}(\omega))\big).
    \ens
    This model assumes we make a single decision now, $\xtilde_{tt}$, then observe a complete sample path of realizations:
    \bns
    (\Wtilde_{t,t+1}(\omega),\Wtilde_{t,t+2}(\omega), \ldots, \Wtilde_{t,t+H}(\omega)),
    \ens
    and then make a complete set of decisions for each sample path $\omega$:
    \bns
    (\xtilde_{t,t+1}(\omega),\xtilde_{t,t+1}(\omega),\ldots, \xtilde_{t,t+H}(\omega)).
    \ens
    This approach insures that the decision now, $x_t=\xtilde_{tt}$, does not depend on what outcome happens in the future, but future decisions, $\xtilde_{tt'}$ for $t' > t$ are allowed to see the entire history of future information.  Next we create a set of scenarios of $\Wtilde_{tt'}(\omega)$ that we denote $\Omegatilde_t$.  Our policy is then written
    \bn
    \hspace{-.15in}X^\pi_t(S_t) = \argmin_{\xtilde_{tt},(\xtilde_{tt'}(\omega))_{t'=t+1}^{t+H},\forall \omega\in\Omega_{[t,t+H]}} \left(\Ctilde(\Stilde_{tt},\xtilde_{tt}) + \hspace{-.05in} \sum_{\omega\in\Omega_[t,t+H]} \Ctilde(\Stilde_{tt'}(\omega),\xtilde_{tt'}(\omega))\right)\hspace{-.04in}. \label{eq:twostagepolicy}
    \en
    The optimization problem in \eqref{eq:twostagepolicy} is typically around $|\Omegatilde_t|$ times bigger than the deterministic problem in equation \eqref{eq:deterministicdla}, but at least it is solvable.
\end{itemize}
The deterministic lookahead, which is also known as ``model predictive control'' in the optimal control literature, is largely dismissed by the stochastic optimization community as little more than a ``deterministic approximation.''  By contrast, the policy based on the two-stage stochastic program in \eqref{eq:twostagepolicy} has appeared in thousands of academic publications, often without recognizing that it is a suboptimal policy for the original optimization problem in equation \eqref{eq:baseobjective}.  See \cite{PowellUnitComm2019} for a discussion of the limitation of scenario trees for the stochastic unit commitment problem.

\subsection{Discussion}
Policies based on solving lookahead models depend on the accuracy of the model to produce good decisions.  The problem is that the solution of full, multistage stochastic decision problems is inherently intractable, forcing the use of very strong approximations such as two-stage stochastic programs.

We propose to extend the idea of a parametric cost function approximation, which we first introduced for state or single-period problems (in section \ref{sec:cfas}), to deterministic lookahead models.  Then, instead of depending on developing an accurate stochastic lookahead model, we exploit structure in the problem to parameterize the deterministic lookahead model to produce behaviors that make the optimal solution more robust.  We then depend on the tuning using a realistic, stochastic simulator (equation \eqref{eq:baseobjective}) to produce the best values of the parameters.

Some advantages of this approach include:
\begin{itemize}
    \item The tuning is done in a realistic simulator (equation \eqref{eq:baseobjective}) that does not need to make simplifications such as an exogenous information process that is independent of decisions.
    \item The simulator can capture any level of detail in the dynamics of the system.
    \item The parameterization of the policy can exploit structure and the modeler's intuition about how uncertainty is likely to affect the solution (an assumption that is made in virtually all parametric models in machine learning).
    \item The resulting policy will generally have the same computational demands as a classical (unparameterized) deterministic lookahead, which would be much easier than solving any stochastic lookahead model.
%    \item Particularly important is that we never have to worry about whether the approach will be used in practice; there are many settings where a parameterized deterministic lookahead is already being used, but may be improved with either a better parameterization, in addition to tuning of the parameters.
    \item This idea has been widely used in industry in an ad-hoc manner.  Specifically, industrial applications will insert parameters without a) recognizing that they are creating a class of policy that is a solution (albeit a suboptimal one) to the optimization problem \eqref{eq:baseobjective} and b) without recognizing that the parameters need to be tuned using the framework of \eqref{eq:baseobjective}.
\end{itemize}
At the same time, we have to consider:
\begin{itemize}
    \item We need to have the intuition into how uncertainty changes the solution we would get from a deterministic lookahead model.
    \item Despite over 60 years of research into stochastic search, parameter tuning remains difficult, but the difficult part is in the research lab where it belongs, not in the field.
\end{itemize}

A major goal of this paper is bring to the attention of the research community in stochastic optimization that a parameterized deterministic lookahead is as valid an approach to the stochastic optimization in \eqref{eq:baseobjective} as any policy based on a stochastic lookahead.  We believe that there are problems where the parameterized deterministic lookahead, in addition to its computational advantages, may outperform a two-stage stochastic program in terms of its ability to solve \eqref{eq:baseobjective}.

\section{Examples of cost function approximations}
\label{sec:examplecfas}
Recall that there are two ways to parameterize an optimization problem: in the objective function, as we did in equation \eqref{eq:linearobjcorrection}, and in the constraints as we did in equations \eqref{eq:constraintcorrection1} - \eqref{eq:constraintcorrection3}.  In this section we are going to provide more concrete examples.

\subsection{CFAs for dynamic assignment problems}

The truckload trucking industry matches drivers to loads over time.  Let
\bns
x_{td\ell} &=& \textwrap{1 if we assign driver $d$ to load $\ell$ at time $t$, 0 otherwise,}\\
c_{td\ell} &=& \textwrap{the contribution of assigning driver $d\in\Dcal_t$ to load $\ell\in\Lcal_t$ at time $t$, including the revenue generated by the load, the cost of moving empty to the load, as well as penalties for late pickup or delivery.}
\ens
We can perform a myopic assignment of drivers to loads by solving
\bn
X^{Assign}(S_t) = \argmax_{x_t} \sum_{d\in\Dcal_t} \sum_{\ell\in\Lcal_t} c_{td\ell} x_{td\ell}. \label{eq:cfadriver}
\en
A potential problem with a myopic policy is that there may be loads that are not assigned and are then held in the hope that a driver may be found to move the load at a later time.  However, the load may be in a location where we do not traditionally have drivers.  We can create an artificial incentive.  Let
\bns
\tau_{t\ell} &=& \textwrap{the amount of time that load $\ell$ has been held at time $t$.}
\ens
Now consider the following modified optimization problem
\bn
X^{Assign}(S_t|\theta) = \argmin_{x_t} \sum_{d\in\Dcal_t} \sum_{\ell\in\Lcal_t} (c_{td\ell} -\theta \tau_{t \ell})x_{td\ell}. \label{eq:cfadriver2}
\en
$X^{Assign}(S_t|\theta)$ is now a parameterized cost function approximation with a modified objective function.

\subsection{A dynamic shortest path problem}
\label{sec:dynamicshortestpath}
Everyone is familiar with the process of navigation systems repeatedly solving shortest path problems to a destination as it receives updates to estimates of travel times around the network.  This is, of course, a fully sequential decision problem that can be modeled as a dynamic program.  We can model the problem as a sequential decision problem.

Let $R_{t}$ be the location of the traveler at time $t$, and let $\cbar_{tij}$ be our estimate of the cost of traversing link $(i,j)$ given what we know at time $t$.  The estimates $\cbar_t$ evolve over time according to
\bn
\cbar_{t+1,ij} = \cbar_{t,ij} + \delta \cbar_{t+1,ij}. \label{eq:shortestpathdynamictransition1}
\en
where $\delta \cbar_{t+1,ij}$ is the change in the estimate of $\cbar_{tij}$ given new observations of traffic.

We assume that at each time $t$ we are at an intersection where we have to make a decision given by
\bns
x_{tij} &=& \left\{\begin{tabular}{cl} 1 & \mbox{if we traverse link $i$ to $j$ when we are at $i$ at time $t$,}\\
                                      0 & \mbox{otherwise.} \end{tabular} \right.
\ens
We make this decision using a policy $X^\pi(S_t)$, where the state $S_t$ is given by
\bns
S_t = (R_t, \cbar_t).
\ens
Models of shortest path problems typically overlook the need to include the vector of estimates $\cbar_t$ in the state variable.  This is precisely why we cannot solve this problem (at least not optimally) using classical shortest path algorithms.

Our challenge is to then find the best policy $X^\pi(S_t)$ that solves
\bn
\min_\pi F^\pi(\theta) = \E \left\{\sum_{t=0}^T \sum_{(i,j)} \chat_{tij}X^\pi(S_t|\theta)|S_0\right\}, \label{eq:shortestpathdynamicobjective}
\en
where $\chat_{tij}$ is the actual cost we experience traversing link $(i,j)$ at time $t+1$.

A natural strategy is to fix the vector of estimates of link costs $\cbar_t$ and solve a shortest path problem to the destination, updating the shortest path as $\cbar_t$ evolves to $\cbar_{t+1}$ (and the traveler transitions to a new node).  This is a lookahead policy based on a lookahead model that uses fixed estimates $\cbar_t$ rather than modeling their stochastic evolution.
%Figure \ref{fig:cfadynamicshortestpath} illustrates the use of a lookahead policy for this problem.
%\begin{figure}[tb]
%\begin{center}
%\begin{tabular}{c}
%\includegraphics[width=5.0in]{cfadynamicshortestpath}
%\end{tabular}
%\caption{Illustration of rolling solution of deterministic shortest path problems using costs $\ctilde^\pi_t(\theta)$.}
%\label{fig:cfadynamicshortestpath}
%\end{center}
%\end{figure}

The question is: Can we do better?  A limitation of the classical approach of solving sequences of deterministic lookaheads is that it fails to recognize that some links can have long tails, which introduces the risk of arriving late.  An alternative is to replace $\cbar_{tij}$ with the $\theta$-percentile of the distribution for each link.  Let
\bns
\cbar^\pi_{t,ij}(\theta) &=& \textwrap{the $\theta$-percentile of the travel time for link $(i,j)$ given our estimate at time $t$.}
\ens
This means we still have a deterministic shortest path problem, but now we have to tune $\theta$ using \eqref{eq:shortestpathdynamicobjective} using our tools from stochastic search.

\section{An energy storage example with rolling forecasts}
\label{sec:cfaenergystorage}
\begin{figure}[tb]
\begin{center}
    \includegraphics[width=4.5in]{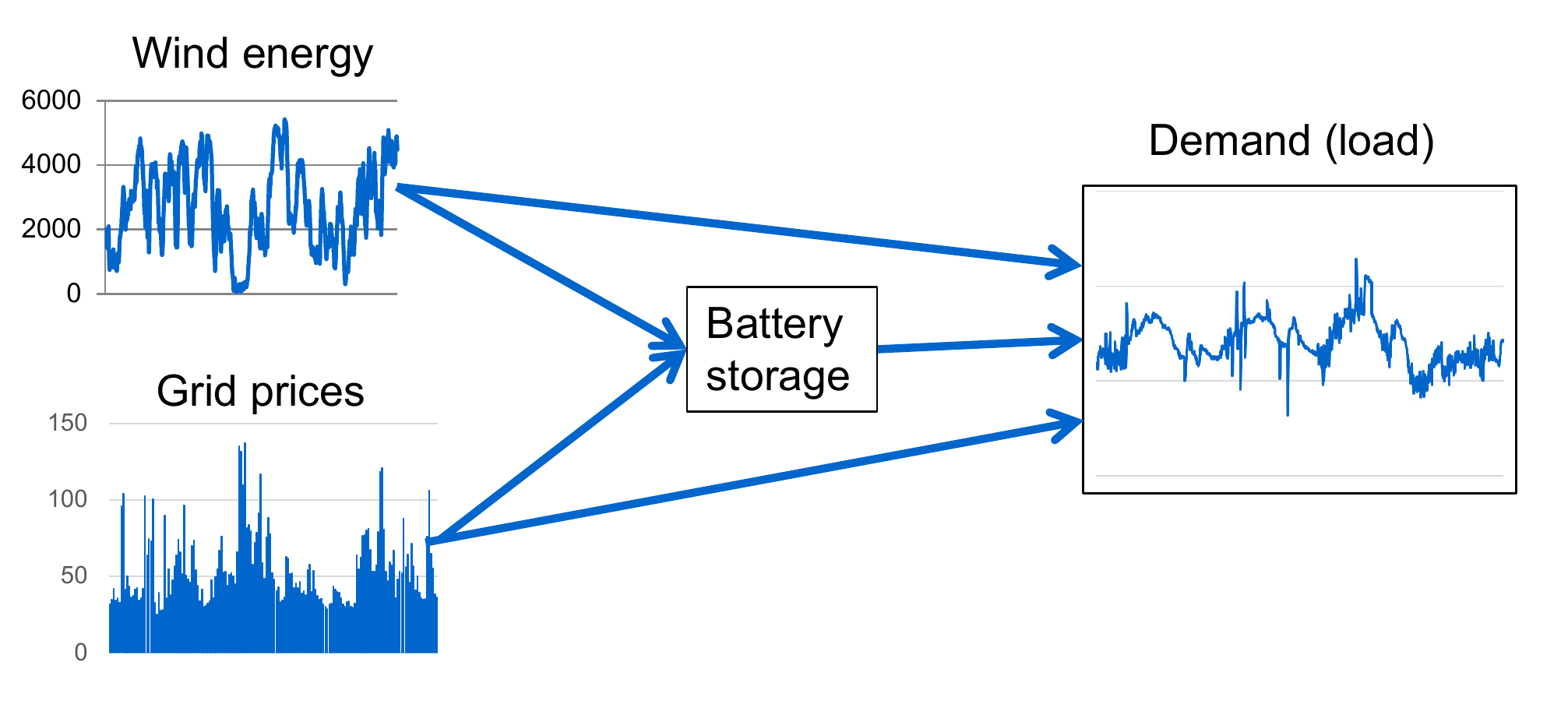}
    \caption{Energy storage system, including a renewable source (wind), energy from the grid at real-time prices, battery storage, and a load.}
    \label{fig:energysystemnew}
\end{center}
\end{figure}

One of the most overlooked modeling issues in operations research is the proper handling of rolling forecasts.  We use the setting of an energy storage system (depicted in figure \ref{fig:energysystemnew}) which draws energy from a wind farm or the power grid to serve a time-varying load, with a finite capacity storage device (and fixed transmission constraints) to help smooth the variations.

\begin{figure}[b]
\center{
    \begin{tabular}{cc}
    \includegraphics[width = 2.70in]{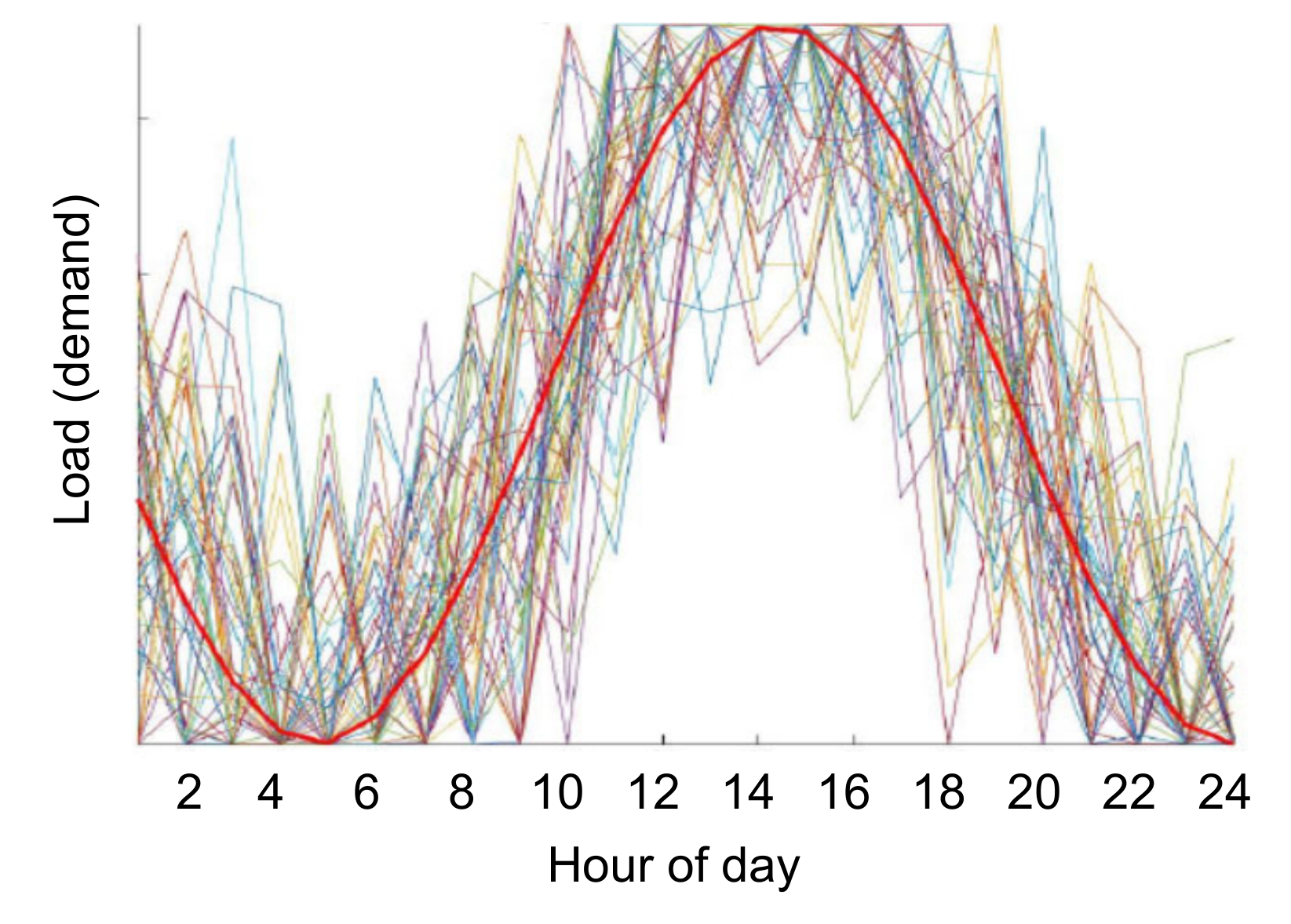} &
    \includegraphics[width = 3.40in]{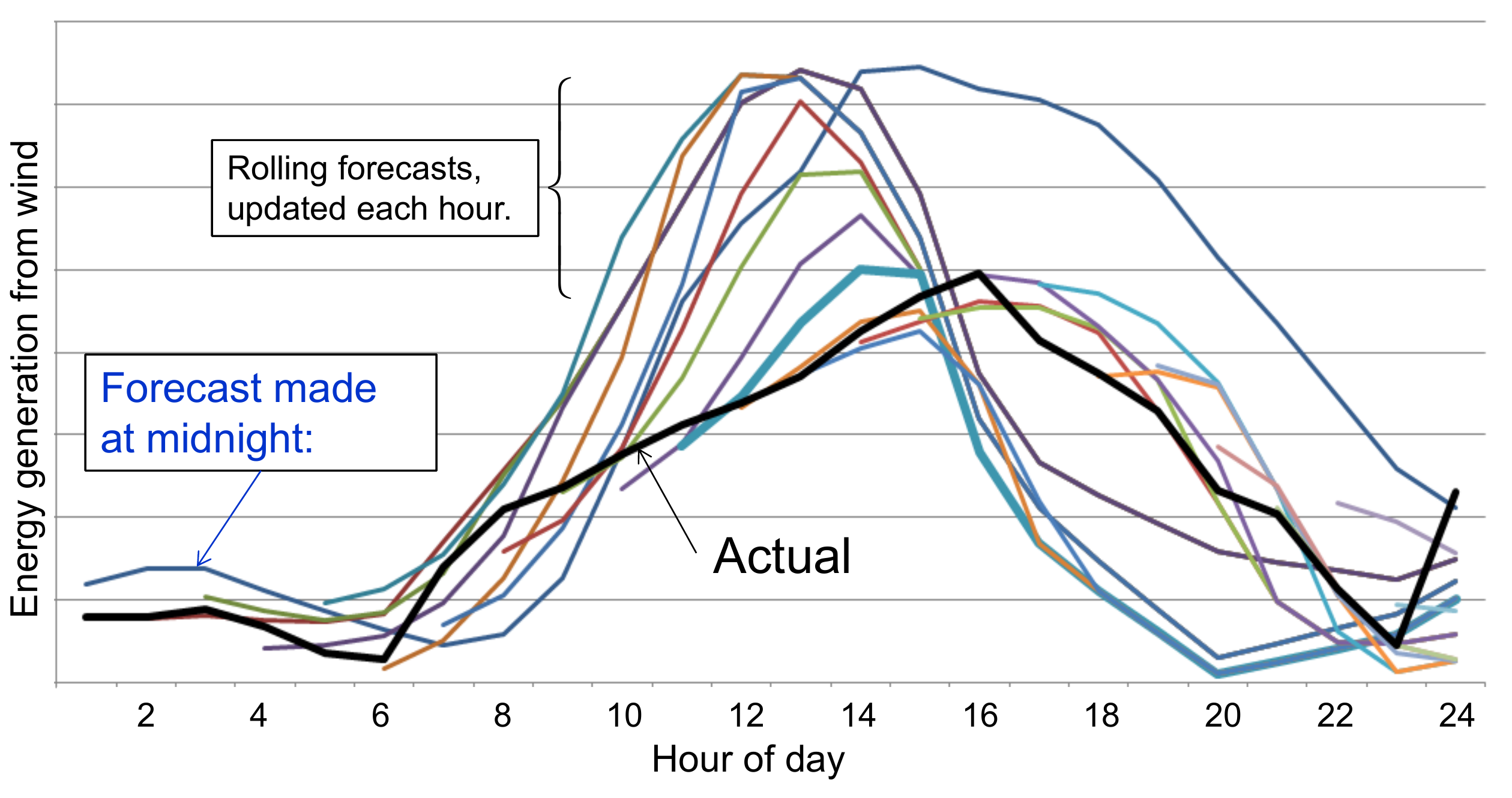} \\
                    (a)                  &                    (b)
    \end{tabular}
    \caption{(a) Energy load by hour of day and (b) rolling forecast, updated hourly.}
    \label{fig:cfahourofday}
    }
\end{figure}

Our energy system has some important characteristics that make it unusually difficult as a stochastic optimization problem:
\begin{itemize}
    \item The energy demands follow a highly time-dependent diurnal pattern (see figure \ref{fig:cfahourofday}(a)).
    \item The energy from wind is highly stochastic.  We have rolling forecasts, updated every hour, but these rolling forecasts evolve considerably over time, as indicated in figure \ref{fig:cfahourofday}(b).  Our rolling wind forecast data was obtained courtesy of PJM Interconnections.
    \item There is unlimited power available from the grid, but at highly stochastic prices.  We can buy from, and sell to, the grid.
    \item The battery has fixed capacity, and the transmission lines are also capacitated, which limits our ability to transmit and store power.  For this reason, the ability to anticipate surges and dips in wind energy requires that we be able to plan into the future.
\end{itemize}

We use a {\it martingale model of forecast evolution} (MMFE) \citep{Heath1994,Graves1986} where forecasts (for energy from the wind farm and the demand) evolve according to
\bn
f^E_{t+1,t'} = f^E_{tt'} + \varepsilon^E_{t+1,t'}, \label{eq:rollingforecastE}\\
f^D_{t+1,t'} = f^D_{tt'} + \varepsilon^D_{t+1,t'}. \label{eq:rollingforecastD}
\en
where $\varepsilon^E_{t+1,t'} \sim N(0,\sigma^2_E)$ and $\varepsilon^D_{t+1,t'} \sim N(0,\sigma^2_D)$ represents the exogenous change in the forecast of energy from the wind farm and the demand for time $t'$.

\subsection{A model of the energy system}

Our model consists of five elements: state variables $S_t$, decision variables $x_t$, exogenous information variables $W_{t+1}$, the transition function $S_{t+1}=S^M(S_t,x_t,W_{t+1})$, and the objective function.  These are given below:

{\noindent \bf State variables:}- The state of the system at time $t$ is all the information we need to model our system from time $t$ onward, which means the information need to compute costs and constraints, making a decision, and compute the transition function.  For the energy problem this information is:
\bns
D_t   &=& \textwrap{Demand (``load'') for power during hour $t$.}\\
E_t   &=& \textwrap{Energy generated from renewables (wind/solar) during hour $t$.}\\
R_t   &=& \textwrap{Amount of energy stored in the battery at time $t$.}\\
u_t   &=& \textwrap{Limit on how much generation can be transmitted at time $t$ (this is known in advance).}\\
p_t   &=& \textwrap{Price to be paid for energy drawn from the grid at time $t$.},\\
f^D_{tt'} &=& \textwrap{Forecast of $D_{t'}$ made at time $t$.}\\
f^E_{tt'} &=& \textwrap{Forecast of $E_{t'}$ made at time $t$.}
\ens
These variables make up our state variable:
\bns
S_t = (R_t, D_t, E_t,  (f^D_{tt'})_{t'\geq t}, (f^E_{tt'})_{t'\geq t}).
\ens
Rolling forecasts are widely used in dynamic models, but the recognition that the forecast itself belongs in the state variable has been recognized by only a small handful of authors, including \cite{Chen1999},\cite{Iida2006} and \cite{Lai2008}.

\noindent {\bf Decision variables:} - These are the flows between each of the elements of our energy system:
\bns
x_{t}      &=& \textwrap{Planned generation of energy during hour $t$ which consists of the following elements:}\\
x^{ED}_{t} &=& \textwrap{flow of energy from wind to demand,}\\
x^{EB}_{t} &=& \textwrap{flow of energy from wind to battery,}\\
x^{GD}_{t} &=& \textwrap{flow of energy from grid to demand,}\\
x^{GB}_{t} &=& \textwrap{flow of energy from grid to battery,}\\
x^{BD}_{t} &=& \textwrap{flow of energy from battery to demand.}
\ens
We would normally write out the constraints that these flows have to satisfy.  These consist of the flow conservation constraints, as well as upper bounds due to transmission constraints, in addition to nonnegativity constraints on all the variables except $x^{GB}_{t}$ since energy is allowed to flow both ways between the grid and the battery.  For compactness, we are going to represent the constraints using
\bns
A_t x_t &=& R_t,\\
x_t     & \leq & u_t,\\
x_t     & \geq & 0.
\ens

\noindent {\bf Exogenous information} - For the variables with forecasts (demand and wind energy), the exogenous information is the change in the forecast, or the deviation between forecast and actual:
\bns
\varepsilon^D_{t+1,\tau} &=& \textwrap{Change in the forecast of demand (for $\tau > 1$ periods in the future) that we first learn at time $t+1$, or the deviation between actual and forecast (for $\tau = 1$).}\\
\varepsilon^E_{t+1,\tau} &=& \textwrap{Change in the forecast of wind energy (for $\tau > 1$ periods in the future) that we first learn at time $t+1$,  or the deviation between actual and forecast (for $\tau = 1$).}
\ens
We assume that prices evolve purely exogenously with deviations:
\bns
\phat_{t+1} &=& \textwrap{Change in grid prices between $t$ and $t+1$.}
\ens
Our exogenous information is then
\bns
W_{t+1} = ((\varepsilon^D_{t+1,\tau},\varepsilon^E_{t+1,\tau})_{\tau \geq 1}, \phat_{t+1}).
\ens

\noindent {\bf Transition function} - The variables that evolve exogenously are
\bns
f^D_{t+1,t'} &=& f^D_{tt'} + \varepsilon^D_{t+1,t'-t-1},~~t'=t+2, \ldots,\\
D_{t+1}      &=& f^D_{t+1,t'} + \varepsilon^D_{t+1,1},\\
f^E_{t+1,t'} &=& f^E_{tt'} + \varepsilon^E_{t+1,t'-t-1},~~t'=t+2, \ldots,\\
E_{t+1}      &=& f^E_{t+1,t'} + \varepsilon^E_{t+1,1},\\
p_{t+1}      &=& p_t + \phat_{t+1}.
\ens
The energy in storage evolves according to
\bns
R_{t+1} = R_t  + x^{EB}_t + x^{GB}_t - x^{BD}_t.
\ens
These equations make up our transition function $S_{t+1} = S^M(S_t,x_t,W_{t+1})$.\\

\noindent {\bf Objective function} - Our single-period contribution function is
\bns
C(S_t,x_t) = p_t \big(x^{GB}_t + x^{GD}_t\big).
\ens
Our objective function, then, would be
\bn
\max_{\pi=(f\in\Fcal,\theta\in\Theta^f)} F^\pi(\theta) = \E\left\{\sum_{t=0}^T C(S_t,X^\pi(S_t|\theta))|S_0\right\}, \label{eq:cfaenergyobjective}
\en
where $S_{t+1} = S^M(S_t,x_t=X^\pi(S_t|\theta),W_{t+1})$, and given a model of the uncertainty that enters our system through the initial state $S_0$ and the exogenous information sequence $W_1, W_2, \ldots, W_T$.  As in the past, we can estimate this objective function by simulating our policy, which we present next.

\subsection{Designing a policy}

%xxx Be sure to discuss reasons why affine policies, VFAs and stochastic lookaheads will not work.\\
%xxx Point out that a full stochastic lookahead model has to model the evolution of forecasts, as well as the evolution of decisions.  Introduce $f^D_{t',t''}$ and $\xtilde_{t',t''}$ in the lookahead model.\\
%xxx Note that there is no pattern among $\theta_\tau$ across $\tau$.

There is a very small literature that addresses inventory problems while explicitly recognizing rolling forecasts $f_t = (f^D_{tt'},f^E_{tt'})_{t'\geq t}$.  \cite{Iida2006} shows that an order-up-to policy parameterized by $\theta(f_t)$ is optimal, but does not attempt to compute the multidimensional function $\theta(f_t)$.  \cite{Lai2008} considers price forecasts in the context of natural gas storage, formulating the dynamic program with $f_t$ in the state variable and showing that an order-up-to policy $\theta(f_t)$ is optimal  but also proposes an approximation based on supporting hyperplanes.  \cite{Chen1999} formulates an inventory problem with rolling forecasts and attempts to use approximate dynamic programming, but is limited to about a half dozen dimensions.  None of these papers considers the much more difficult problem of bounds on order quantities and storage, which makes the problem much harder, and invalidates the optimality proofs of order-up-to policies that can be written as $\theta(f_t)$ (they have to be time-dependent, and state-dependent $\theta_t(S_t)$ where $S_t$ includes all the elements of the state variable).

Drawing on the framework of a parametric cost function approximation, we are going to start with a classical deterministic lookahead model.  We begin by creating the decision variables for our lookahead model
\bns
\xtilde_{tt'} = (\xtilde^{ED}_{tt'},\xtilde^{EB}_{tt'},\xtilde^{GD}_{tt'},\xtilde^{GB}_{tt'},\xtilde^{BD}_{tt'}), ~t+1 \leq  t' \leq t+H,
\ens
which parallels the elements of $x_t$ for each time $t'$ in the future.

This is a time-dependent problem with complex interactions between the uncertain supply of wind energy, the price of energy from the grid, and the time-dependent nature of demands that have to be satisfied over a capacitated grid.  It seems natural to start by creating a policy based on a deterministic lookahead model given by:

{\small
\bn
X^{DLA}(S_t) \hspace{-.10in}& = &\hspace{-.10in} \argmax_{x_t,(\xtilde_{tt'},t'=t+1, \ldots, t+H)} \left(p_t (x^{GB}_t + x^{GD}_t) + \sum_{t'=t+1}^{t+H} \ptilde_{tt'} (\xtilde^{GB}_{tt'} + \xtilde^{GD}_{tt'}) \right) \nonumber \\
& &\label{eq:DLA0}
\en
}
subject to the following constraints: First, for time $t$ we have
\bn
x^{BD}_{t} - x^{GB}_{t} - x^{EB}_{t} & \leq & R_{t}, \label{eq:DLA1} \\
\Rtilde_{t,t+1}- (x^{GB}_{t} + x^{EB}_{t} - x^{BD}_{t}) &=& R_{t}, \label{eq:DLA2}\\
x^{ED}_{t} + x^{BD}_{t} + x^{GD}_{t} & = & D_{t},\label{eq:DLA3}\\
x^{EB}_{t} + x^{ED}_{t}  & \leq & E_{t}, \label{eq:DLA4} \\
x^{GD}_t, x^{EB}_t, x^{ED}_t, x^{BD}_t & \geq & 0, \label{eq:DLA5}
%\Rtilde_{tt} &=& R_t - (x^{BD}_{t} - x^{GB}_{t} - x^{EB}_{t}). \label{eq:DLA6}
\en
Then, for $t' = t+1, \ldots, t+H$ we have
\bn
%\Rtilde_{t,t'+1} &=& \Rtilde_{t,t'} - (\xtilde^{BD}_{tt'} - \xtilde^{GB}_{tt'} - \xtilde^{EB}_{tt'}),\label{eq:DLA6}\\
\xtilde^{BD}_{tt'} - \xtilde^{GB}_{tt'} - \xtilde^{EB}_{tt'} & \leq & \Rtilde_{tt'}, \label{eq:DLA6} \\
\Rtilde_{t,t'+1}- (\xtilde^{GB}_{tt'} + \xtilde^{EB}_{tt'} - \xtilde^{BD}_{tt'}) &=& \Rtilde_{tt'}, \label{eq:DLA7}\\
\xtilde^{ED}_{tt'} + \xtilde^{BD}_{tt'} + \xtilde^{GD}_{tt'} & = & f^D_{tt'},\label{eq:DLA8}\\
\xtilde^{EB}_{tt'} + \xtilde^{ED}_{tt'}  & \leq & f^E_{tt'}. \label{eq:DLA9}
\en

The weakness in this model is the forecasts of wind energy $f^E_{tt'}$ and demand $f^D_{tt'}$.  One idea would be to ``discount'' these forecasts by multiplying each forecast with coefficients $\theta^E_{t'-t}$ and $\theta^D_{t'-t}$ that depend on how far into the future we are trying to forecast.  Using this approach, we create a parameterized policy by replacing equations \eqref{eq:DLA8} and \eqref{eq:DLA9} with
\bn
\xtilde^{ED}_{tt'} + \xtilde^{BD}_{tt'} + \xtilde^{GD}_{tt'} & = & \theta^D_{t'-t}f^D_{tt'},\label{eq:DLA8a}\\
\xtilde^{EB}_{tt'} + \xtilde^{ED}_{tt'}  & \leq & \theta^E_{t'-t}f^E_{tt'}. \label{eq:DLA9a}
\en

%It is easy to criticize this approach as just a heuristic, deterministic lookahead.  This attitude completely ignores a) the intuition behind the parameterization and b) the fact that we are tuning $\theta$ (using equation \eqref{eq:cfaenergyobjective}) with a stochastic, dynamic model that does not make any of the approximations that we would need to create a solvable stochastic lookahead model.

\begin{figure}[tb]
\begin{center}
\begin{tabular}{c}
\includegraphics[width=3.5in]{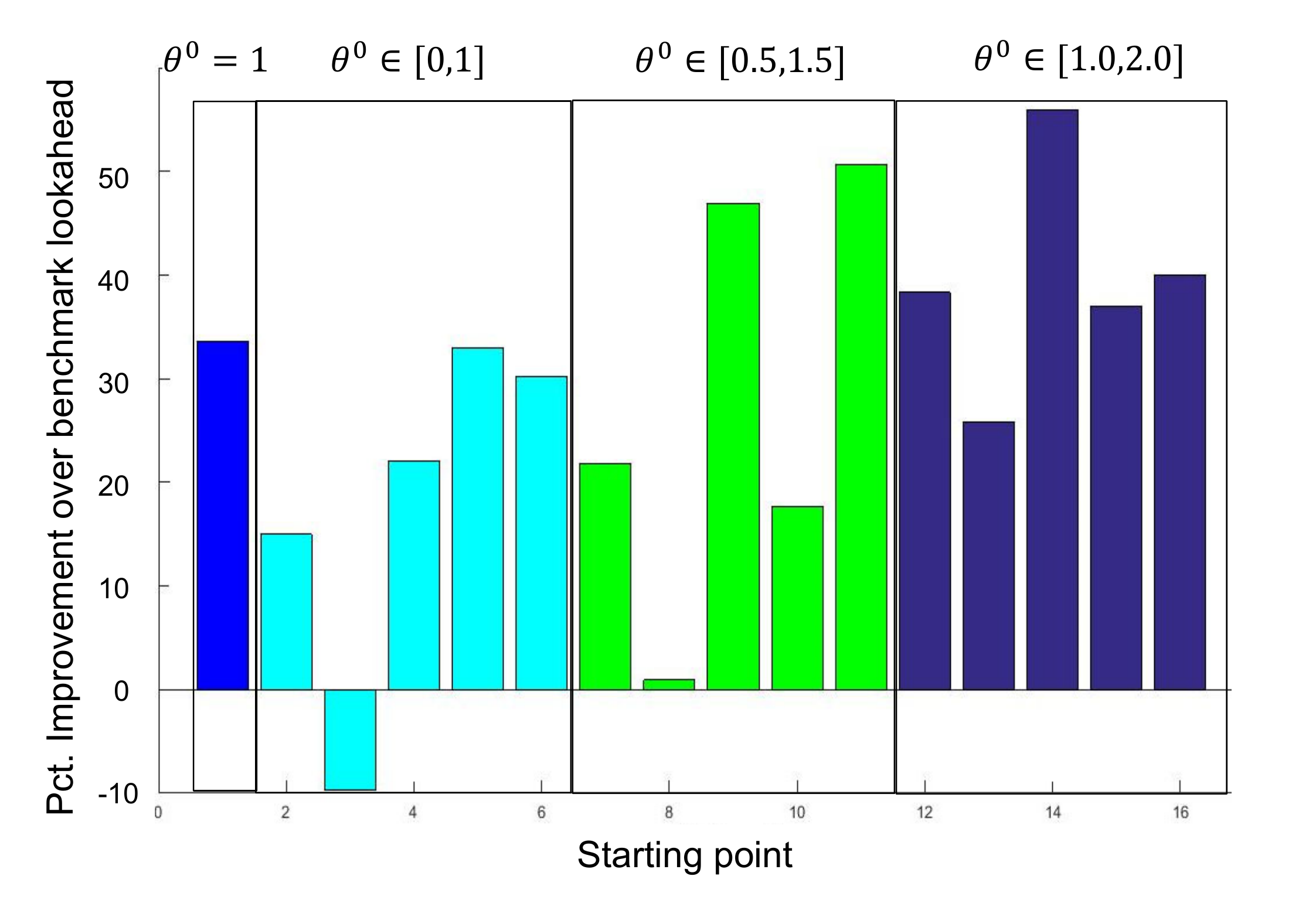}
\end{tabular}
\caption{Relative performance of optimized parameterized lookahead policy to the performance of a deterministic lookahead with $\theta = 1$.}
\label{fig:cfaenergysystemperformance}
\end{center}
\end{figure}

The problem of tuning the parameter vector $\theta$ is not easy, but there are a number of strategies we can draw on.  We designed an algorithm \citep{Ghadimi2022} using a stochastic gradient algorithm based on Spall's simultaneous perturbation stochastic approximation algorithm \cite{spall2005introduction} which is well suited to problems with multidimensional parameters.  We then compared the performance of the optimized parametric policy against a base policy with $\theta = 1$.  The results are shown in figure \ref{fig:cfaenergysystemperformance} for optimized $\theta$ using a range of different starting points, where the left-most bar uses an initial starting point of $\theta^0 = 1$.  Most of the results show improvements of 20 to 50 percent.

We claim that the roughly 30 percent improvement is quite significant, given that it does not come at any additional computational cost in the field.  At the same time we observe that there is no alternative computational strategy that would be guaranteed to be better.  Two-stage stochastic programs do not even attempt to model the evolution of estimates of the forecasts.  Approximate dynamic programming would not be able to capture the complex nonlinearities of forecasts in the state variable, especially for the capacitated problem that we are solving.

\section{Closing remarks}
\label{sec:closingremarks}
The major goal of this paper is to make the argument that a parameterized deterministic optimization model is a perfectly valid basis for a policy for (stochastic) sequential decision problems.  The research community needs to accept that tuning a parameterized policy using a stochastic base model such as that given in equation \eqref{eq:baseobjective} is a form of stochastic optimization, even if the policy requires solving a deterministic optimization problem.

Deterministic lookahead models are widely used in practice because they easily handle complexity, and are relatively easy to solve.  The use of parameterized cost function approximations enjoys several significant strengths, especially in the context of complex problems.  Some examples are:
\begin{itemize}
\item The parametric cost function approximations naturally handles the dynamics of a highly time-dependent problem with complicating constraints (as we encountered in the energy storage problem) that emphasize the importance for planning into the future.
\item Although a time-dependent problem requires time-dependent behavior, the effect of incorporating a rolling forecast produces a stationary policy.  The function $X^\pi(S_t)$ is not time dependent, and the vector $\theta$ is not time-dependent.  This property significantly simplifies the search process for $\theta$.
%\item The policy naturally captures the dependence on all the elements of the state variable $S_t$, eliminating the need to design a complex function such as a state-dependent order-up-to parameter $\theta_t(S_t)$.
\item Rolling forecasts arise in many settings, yet introduce complex stochastic interactions between forecasts and decisions in the future, which impact the decisions made now.  Capturing this in a stochastic lookahead policy is exceptionally difficult, but is quite easy in a simulator.
\item Parametric deterministic lookahead policies can capture complex state variables (the rolling forecast is just one example) much more easily than policies based on stochastic lookaheads.  Similarly, simulating complex state variables for the purpose of parameter tuning is also quite easy.
\item The parametric deterministic lookahead policy is generally easy to compute in the field.
\end{itemize}

While this approach is both attractive (since it is easy to implement) and promising (see figure \ref{fig:cfaenergysystemperformance}), serious research issues remain:
\begin{itemize}
  \item Designing the best parameterization is difficult, but closely parallels the challenges of model design in machine learning.  For our energy storage problem, we might say that multiplying coefficients $\theta$ times the forecasts is intuitively appealing, but other parameterizations are possible, such as ensuring that the energy in storage in future time periods stays above a minimum level (as a reserve) and below a maximum level (so we can store unexpected surges in wind).  We suspect that most industrial applications at best use intuitive parameterizations without the benefits of experimental testing in a simulator.
%  \item Even if we accept the intuition of ``discounting'' forecasts, it seems reasonable to expect that $\theta_\tau$ (for either forecasted wind energy or demand) would trace a smooth function in $\tau$.  One set of optimized values for $\theta^E_\tau$ is shown in figure \ref{fig:optimizedthetas}, which does not show any pattern at all.  It would be interesting to understand if this is a result of noise in the stochastic search algorithm, or if it reflects the properties of the problem.
  \item Stochastic search remains a challenge. For example, the objective function \eqref{eq:cfaenergyobjective} is nonconvex in $\theta$ (when we limit the search over $\theta$).  Simulating policies can also be quite noisy.
  \item While it is natural to assume that we can tune the parameters using a simulator, there will be many settings where a simulator is not available (or would not be trusted).  An important research challenge is to perform online parameter tuning in the field so that the policy adapts to changing conditions.
\end{itemize}

We hope that the thoughts in this paper encourage the stochastic optimization community to include parameterized deterministic models as valid policies for stochastic optimization problems.  This initiative is likely to be warmly endorsed by industry that is already implementing parameterized deterministic models, without the benefits of careful design of the parameterization and parameter tuning.

%{\bf Acknowledgments} \\ \\
%This research has been supported in part by the National Science Foundation, grant CMMI-1537427.

%\singlespace
%\bibliographystyle{../../agsm}
%\bibliography{ghadimilibrary,../../bib2/library}

%\newpage

%\input{Appendix}

\end{document}